\documentclass[smallextended]{svjour3}       
\smartqed  
\usepackage{graphicx}
\usepackage[cmex10]{amsmath}
\usepackage{amssymb}
\usepackage{amsbsy}
\usepackage{algorithmicx,algorithm}
\usepackage{algpseudocode}
\usepackage{threeparttable}
\usepackage{enumerate}
 \usepackage{mathptmx}      
\newtheorem{thm1}{Theorem}
\newtheorem{thm2}[thm1]{Theorem}
\newtheorem{thm3}[thm1]{Theorem}
\newtheorem{thm4}[thm1]{Theorem}
\newtheorem{thm5}[thm1]{Theorem}
\newtheorem{thm6}[thm1]{Theorem}

\newtheorem{lma1}{Lemma}
\newtheorem{lma2}[lma1]{Lemma}
\newtheorem{lma3}[lma1]{Lemma}
\newtheorem{lma4}[lma1]{Lemma}

\newtheorem{cor1}{Corollary}
\newtheorem{cor2}[cor1]{Corollary}
\newtheorem{cor3}[cor1]{Corollary}
\newtheorem{cor4}[cor1]{Corollary}
\newtheorem{cor5}[cor1]{Corollary}

\begin{document}

\title{New Kloosterman sum identities from the Helleseth-Zinoviev result on $ Z_{4}$-linear Goethals codes 
\thanks{Communicated by Pascale Charpin, Alexander Pott, Dieter Jungnickel}
}


\author{Minglong Qi \and Shengwu Xiong  }


\institute{Minglong Qi \and Shengwu Xiong  \at
              School of Computer Science and Technology,  Wuhan University of Technology \\
              Mafangshan West Campus, 430070 Wuhan City, China\\
              \email{mlqiecully@163.com (Minglong Qi) }\\
              \email{xiongsw@whut.edu.cn (Shengwu Xiong)}\\ 
}

\date{Received: date / Accepted: date}

\maketitle

\begin{abstract}
In the paper of Tor Helleseth and Victor Zinoviev (Designs, Codes and Cryptography, \textbf{17}, 269-288(1999)), the number of solutions of the system of equations from $ Z_{4} $-linear Goethals codes $ G_{4} $ was determined and stated in Theorem 4.  We found that Theorem 4 is wrong for $ m $ even. In this note, we complete Theorem 4, and present a series of new Kloosterman sum identities  deduced from Theorem 4. Moreover, we show that several previously established formulas on the Kloosterman sum identities can be rediscovered from Theorem 4 with much simpler proofs.

\end{abstract}

\section{Introduction}
Let $ m $ be a positive integer,  $q=2^{m}, F:=\mathbb{F}_{q} $ be the finite field of $ q $ elements, $ F^{*}:=F\setminus \lbrace 0\rbrace $,  and $ F^{**}:=F\setminus \lbrace 0, 1\rbrace $. The well known Kloosterman sums \cite{carlitz,charpin} are defined by 
\begin{equation}\label{kloosterman_sum}
K(a)=\sum\limits_{x\in F^{*}}(-1)^{\mathrm{Tr}(ax+1/x)}
\end{equation}
where $ a\in F^{*} $, and $ \mathrm{Tr}(\cdot) $ is the trace function of $ F $ over $ \mathbb{F}_{2} $.

Let $ b,c\in F $. The problem of finding the coset weight distribution of $ Z_{4} $-linear Goethals codes $ G_{4} $,  is transformed into solving the following nonlinear system of equations over $ F $ \cite{tor-victor}:
\begin{equation}\label{G4_system}
\begin{cases}
x+y+z+u &=1\\
u^{2}+xy+xz+xu+yz+yu+zu &=b^{2}\\
x^{3}+y^{3}+z^{3}+u^{3} &=c
\end{cases}
\end{equation}
where $ x, y, z $ and $ u $ are pairwise distinct elements of $ F $. The number of solutions of (\ref{G4_system}) (see p.284 
of \cite{tor-victor}), denoted by $ \mu_{2}(b,c) $, is given by Theorem 4 of \cite{tor-victor}.

We found that Theorem 4 of \cite{tor-victor}, which gives the explicit evaluation of $ \mu_{2}(b,c) $, is wrong for  $ m $ even. It is obvious that the authors of \cite{tor-victor} forgot to take account of the fact, that $ \mathrm{Tr}(1)=1 $ for $ m $ odd and 0 for $ m $ even, in the last step of the proof of Theorem 4. The following is the correct version of Theorem 4 of \cite{tor-victor}:

\begin{thm1}\label{thm1}
Let $ \mu_{2}(b,c) $ be the number of different 4-tuples $ (x, y, z, u) $, where $ x, y, z, u $ are pairwise distinct elements of $ F $, which are solutions to the system (\ref{G4_system}) over $ F $, where $ b, c $ are arbitrary elements of $ F $. 
\begin{enumerate}[(1)]
\item If $ m $ is odd and $ \mathrm{Tr}(c)=1 $ or $ m $ is even and $ \mathrm{Tr}(c)=0 $, then
\begin{equation}\label{eq_nsol01}
\mu_{2}(b,c)=\frac{1}{6}\bigl(q-8+(-1)^{\mathrm{Tr}(b)}(K(k_{1}k_{2})-3)\bigr).
\end{equation}
\item If $ m $ is odd and $ \mathrm{Tr}(c)=0 $ or $ m $ is even and $ \mathrm{Tr}(c)=1 $, then
\begin{equation}\label{eq_nsol02}
\mu_{2}(b,c)=\frac{1}{6}\bigl(q-2-(-1)^{\mathrm{Tr}(b)}(K(k_{1}k_{2})+3)\bigr).
\end{equation}
Where $  k_{1}=b^{2}+c+1$ and  $ k_{2}=b^{2}+b+c+\sqrt{c} $.
\end{enumerate}
\end{thm1}

In this note, we will show that Theorem 4 of \cite{tor-victor} implies not only a series of new Kloosterman sum identities, but  almost all previously discovered formulas on the identities for Kloosterman sums, thanks to the following theorem which is a  corollary of that theorem, too:

\begin{thm2}\label{thm2}
Let $ b, c\in F,  k_{1}=b^{2}+c+1 $ and $ k_{2}=b^{2}+b+c+\sqrt{c} $. If $ k_{1}k_{2}\ne 0 $, then 
\begin{equation}\label{eq_KS_k1k2}
K(k_{1}k_{2})=K(k_{1}k_{2}+k_{2}).
\end{equation}
\end{thm2}

The idea is to affect a special value to $ k_{1} $ or $ k_{2} $ and then substitute them into (\ref{eq_KS_k1k2}) to arrive at a special Kloosterman sum identity. Let $ b, c\in F$,  and $ n, k $ be rational numbers such that $ b^{n}, b^{k}, c^{n}, c^{k} $ are elements of $ F $. In the sequel, $ k_{1} $ and $ k_{2} $ are affected to the following values:
\begin{equation}\label{k12_values}
\begin{cases}
k_{1} &\in \lbrace  cb^{n}, b^{2}c^{n}, cb^{n}+b^{k}, b^{2}c^{n}+c^{k}\rbrace,\\
k_{2} &\in \lbrace cb^{n}+\sqrt{c}, \sqrt{c}b^{n}+c, b^{2}c^{n}+b, bc^{n}+b^{2}\rbrace.
\end{cases}
\end{equation}

The previously discovered formulas on the identities for Kloosterman sums are stated in the following theorem:

\begin{thm3}\label{thm3}
Let $ a\in F^{**} $. Then
\begin{enumerate}[(1)]
\item $ K\bigl(a^{3}(a+1)\bigr)=K\bigl(a(a+1)^{3}\bigr)$ \cite{shin-kumar-helleseth,shin-sung,tor-victor-ksi1-dm,tor-victor-ksi12-ff} (Helleseth-Zinoviev Formula I).
\item $ K\bigl(a^{5}(a+1)\bigr)=K\bigl(a(a+1)^{5}\bigr) $ \cite{tor-victor-ksi12-ff} (Helleseth-Zinoviev Formula II).
\item $ K\bigl(a^{8}(a^{4}+a)\bigr)=K\bigl((a+1)^{8}(a^{4}+a)\bigr) $\cite{hollmann-xiang} (Hollmann-Xiang Formula).
\item $ K\bigl(a/(a+1)^{4}\bigr)=K\bigl(a^{3}/(a+1)^{4}\bigr) $ \cite{shin-kumar-helleseth} (Shin-Kumar-Helleseth Formula).
\end{enumerate}
\end{thm3}

Note that Shin-Kumar-Helleseth Formula can be deduced from Helleseth-Zinoviev Formula I by the variable change $ a=b/(b+1) $. In this note, we will show that except for Helleseth-Zinoviev Formula II, every formula of Theorem \ref{thm3} can be obtained from Theorem 4 of \cite{tor-victor} with help of Theorem \ref{thm2}.

In Section 2, we at first prove Theorem \ref{thm2} with help of Theorem \ref{thm1} (Theorem 4 of \cite{tor-victor}) and two lemmas of \cite{tor-victor-ksi1-dm}, treat in  details some cases of (\ref{k12_values}) for $ k_{1} $ and $ k_{2} $ which bring us new Kloosterman sum identities, prove Helleseth-Zinoviev Formula I and Hollmann-Xiang Formula, and finally generalize Shin-Kumar-Helleseth Formula.

\section{New identities for Kloosterman sums}

\subsection{Proof of Theorem \ref{thm2}}

Before proving Theorem \ref{thm2}, we need several  preliminary lemmas:
 
\begin{lma1}[\cite{carlitz}]\label{lemma1}
Let $ a\in F^{*} $. Then, $ K(a) =K(a^{2})$.
\end{lma1}

\begin{lma2}[\cite{lidl-nieder}]\label{lemma2}
Let $ a, b\in F $. Then,
\begin{enumerate}[(1)]
\item $ \mathrm{Tr}(a)=\mathrm{Tr}(a^{2}),  \mathrm{Tr}(a+b)=\mathrm{Tr}(a)+\mathrm{Tr}(b)$.
\item $ \mathrm{Tr}(1)=1 $ if $ m $ is odd, and $ \mathrm{Tr}(1)=0 $ if $ m $ is even.
\end{enumerate}
\end{lma2}

The symmetry of the solutions of the nonlinear system (\ref{G4_system}) is characterized by the following two lemmas:

\begin{lma3}\cite[Lemma 10]{tor-victor-ksi1-dm}\label{lamma3}
Let $ b, c $ be any elements of $ F $, where $ F $ has order $ 2^{m} $. Let $ m $ be even. Let $ \mu_{2}(b,c) $ denote the number of solutions to system (\ref{G4_system}). Then the following symmetry conditions are valid:
\begin{equation}\label{eq_symmetry_bp1_meven}
\mu_{2}(b,c)=\mu_{2}(b+1,c)
\end{equation}
and 
\begin{equation}\label{eq_symmetry_cp1_meven}
\mu_{2}(b,c)=\mu_{2}(b,c+1).
\end{equation}
\end{lma3}

\begin{lma4}\cite[Lemma 11]{tor-victor-ksi1-dm}\label{lamma4}
Let $ b, c $ be any elements of $ F $, where $ F $ has order $ 2^{m} $. Let $ m $ be odd. Let $ \mu_{2}(b,c) $ denote the number of solutions to system (\ref{G4_system}). Then the following symmetry conditions are valid:
\begin{equation}\label{eq_symmetry_bp1_modd}
\mu_{2}(b,c)+\mu_{2}(b+1,c)=\dfrac{1}{3}\times
\begin{cases}
(2^{m}-2)&\qquad\text{if}\ \mathrm{Tr}(c)=0,\\
(2^{m}-8)&\qquad\text{if}\ \mathrm{Tr}(c)=1,
\end{cases}
\end{equation}
and 
\begin{equation}\label{eq_symmetry_cp1_modd}
\mu_{2}(b,c)+\mu_{2}(b,c+1)=\dfrac{1}{3}\times
\begin{cases}
(2^{m}-2)&\qquad\text{if}\ \mathrm{Tr}(c)=1,\\
(2^{m}-8)&\qquad\text{if}\ \mathrm{Tr}(c)=0.
\end{cases}
\end{equation}
\end{lma4}

Now, we are ready to prove Theorem \ref{thm2}. The idea of proof is to substitute the formulas for $ \mu_{2}(b,c) $ from Theorem \ref{thm1} into the formulas of Lemma \ref{lamma3} for $ m $ even, and into the formulas of Lemma \ref{lamma4} for $ m $ odd.
\begin{proof} [of Theorem \ref{thm2}]
\begin{enumerate}[(I)]
\item Case for $ m $ even.
\begin{enumerate}
\item Subcase that $ \mathrm{Tr}(c)=0 $. 

Note that $ \mathrm{Tr}(1)=0 $ for $ m $ even by Lemma \ref{lemma2}, and $ \mathrm{Tr}(b+1)=\mathrm{Tr}(b)+\mathrm{Tr}(1)=\mathrm{Tr}(b) $. Let $ k_{1}^{'}:=(b+1)^{2}+c+1, k_{2}^{'}:= (b+1)^{2}+(b+1)+c+\sqrt{c}$. It is clear that $ k_{1}^{'}=k_{1}+1, k_{2}^{'}=k_{2} $ and $ k_{1}^{'}k_{2}^{'}=k_{1}k_{2}+k_{2} $, where $ k_{1}=b^{2}+c+1 $ and $ k_{2}=b^{2}+b+c+\sqrt{c} $.
From (\ref{eq_nsol01}) of Theorem \ref{thm1}, we have 
\begin{equation*}
\begin{split}
\mu_{2}(b+1,c)&=\frac{1}{6}\bigl(q-8+(-1)^{\mathrm{Tr}(b+1)}(K(k_{1}^{'}k_{2}^{'})-3)\bigr)\\
&=\frac{1}{6}\bigl(q-8+(-1)^{\mathrm{Tr}(b)}(K(k_{1}k_{2}+k_{2})-3)\bigr)\ \text{and},\\
\mu_{2}(b,c)&=\frac{1}{6}\bigl(q-8+(-1)^{\mathrm{Tr}(b)}(K(k_{1}k_{2})-3)\bigr).
\end{split}
\end{equation*}
Theorem \ref{thm2} follows from (\ref{eq_symmetry_bp1_meven}) of Lemma \ref{lamma3} for this subcase.  

Let $ k_{1}^{''}:=b^{2}+(c+1)+1, k_{2}^{''}:= b^{2}+b+(c+1)+\sqrt{c+1}$. It is clear that $ k_{1}^{''}=k_{1}+1, k_{2}^{''}=k_{2} $ and $ k_{1}^{''}k_{2}^{''}=k_{1}k_{2}+k_{2} $. Since $ \mathrm{Tr}(c+1)=\mathrm{Tr}(c)=0 $, from (\ref{eq_nsol01}) of Theorem \ref{thm1}, we obtain 
\begin{equation*}
\begin{split}
\mu_{2}(b,c+1)&=\frac{1}{6}\bigl(q-8+(-1)^{\mathrm{Tr}(b)}(K(k_{1}^{''}k_{2}^{''})-3)\bigr)\\
&=\frac{1}{6}\bigl(q-8+(-1)^{\mathrm{Tr}(b)}(K(k_{1}k_{2}+k_{2})-3)\bigr)\ \text{and},\\
\mu_{2}(b,c)&=\frac{1}{6}\bigl(q-8+(-1)^{\mathrm{Tr}(b)}(K(k_{1}k_{2})-3)\bigr).
\end{split}
\end{equation*}
Again, Theorem \ref{thm2} follows from (\ref{eq_symmetry_cp1_meven}) of Lemma \ref{lamma3} for this subcase.
\item Subcase that $ \mathrm{Tr}(c)=1 $. 

For this subcase, we use (\ref{eq_nsol02}) of Theorem \ref{thm1} and (\ref{eq_symmetry_bp1_meven}) or (\ref{eq_symmetry_cp1_meven}) of Lemma \ref{lamma3} to prove Theorem \ref{thm2}, and omit the details due to limited space.
\end{enumerate}
\item Case for $ m $ odd. 

For this case, we use (\ref{eq_symmetry_bp1_modd}) or (\ref{eq_symmetry_cp1_modd}) of Lemma \ref{lamma4} and Theorem \ref{thm1} to prove Theorem \ref{thm2}, which is similar to the case for $ m $ even and omitted.
\end{enumerate}
\end{proof}\qed

\subsection{New Kloosterman sum identities and the proof of Theorem \ref{thm3} }


\begin{thm4}\label{thm4-b2cn}
Let $ c\in F $, and  $ n $ be a rational number such that $ c^{n} $ is an element of $ F $. If $ c^{4n}+1\ne 0 $, then 
\begin{equation}\label{eq-thm4-b2cn}
K\biggl(\dfrac{(c^{n+1}+c^{n})^{3}(c^{n+1}+1)}{c^{4n}+1}\biggr)=K\biggl(\dfrac{(c^{n+1}+c^{n})(c^{n+1}+1)^{3}}{c^{4n}+1}\biggr).
\end{equation}
\end{thm4}
\begin{proof}
Let $ k_{1}=b^{2}+c+1=b^{2}c^{n} $. Then, $ b^{2}=(c+1)/(c^{n}+1) $, and $ k_{1}=(c^{n+1}+c^{n})/(c^{n}+1) $. Further, we can obtain 
\begin{equation*}
\begin{split}
k_{2}&=b^2+b+c+\sqrt{c}=\dfrac{c+1}{c^{n}+1}+\sqrt{\dfrac{c+1}{c^{n}+1}}+c+\sqrt{c};\\
k_{2}^{2}&=\dfrac{c^2+1}{c^{2n}+1}+\dfrac{c+1}{c^{n}+1}+c^2+c=\dfrac{(c^{n+1}+c^{n})(c^{n+1}+1)}{c^{2n}+1};\\
(k_{1}k_{2})^{2}&=\dfrac{(c^{n+1}+c^{n})^{3}(c^{n+1}+1)}{c^{4n}+1};
(k_{1}k_{2})^{2}+k_{2}^{2}=\dfrac{(c^{n+1}+c^{n})(c^{n+1}+1)^{3}}{c^{4n}+1}.
\end{split}
\end{equation*}
The actual theorem follows from Lemma \ref{lemma1} and Theorem \ref{thm2}.
\end{proof}\qed

Let $ L(n):= (c^{n+1}+c^{n})^{3}(c^{n+1}+1)/(c^{4n}+1), R(n):= (c^{n+1}+c^{n})(c^{n+1}+1)^{3}/(c^{4n}+1)$, then, $ K(L(n)) =K(R(n))$. We now prove  Helleseth-Zinoviev Formula I by Theorem \ref{thm4-b2cn}.

\begin{proof}[of  Helleseth-Zinoviev Formula I]
\begin{equation*}
\begin{split}
L(1)&=\dfrac{(c^2+c)^3(c^2+1)}{c^4+1}=\dfrac{c^3(c+1)^5}{(c+1)^4}=c^3(c+1),\\
R(1)&=\dfrac{(c^2+c)(c^2+1)^3}{c^4+1}=\dfrac{c(c+1)^7}{(c+1)^4}=c(c+1)^3.
\end{split}
\end{equation*}
Helleseth-Zinoviev Formula I follows from the case $ K(L(1)) =K(R(1)) $ of Theorem \ref{thm4-b2cn}.
\end{proof}\qed


\begin{cor1}\label{cor1-from-thm-b2cn}
Let $ c\in F $ such that the rational functions in $ c $ occurring in the following formulas are valid. Then,
\begin{enumerate}[(1)] 
\item $ K\bigl(c^6(c^2+c+1)/(c+1)^4\bigr)=K\bigl(c^2(c^2+c+1)^3/(c+1)^4\bigr) $,
\item $ K\bigl(c^9(c+1)^3/(c^8+c^4+1)\bigr)=K\bigl(c^3(c+1)^9/(c^8+c^4+1)\bigr) $,
\item $ K\bigl((c+1)^8(c^4+c)\bigr)=K\bigl(c^3(c^3+1)^3\bigr) $,
\item $ K\bigl((c^{11}+c^3)(c^5+1)\bigr)=K\bigl(c(c^5+1)^3\bigr) $,
\item $ K\bigl((c+1)^{20}(c^8+c)\bigr)=K\bigl((c+1)^{4}(c^8+c)^3\bigr) $.
\end{enumerate}
\end{cor1}
\begin{proof}
The first formula of Corollary \ref{cor1-from-thm-b2cn} corresponds to the case $ K(L(2)) =K(R(2))$ of Theorem \ref{thm4-b2cn}, and the second formula to $ K(L(3)) =K(R(3))$, which the proofs are analogous to that of  Helleseth-Zinoviev Formula I and omitted.  We now prove the third formula:
\begin{equation*}
\begin{split}
L(-1/4)&=\dfrac{(c^{3/4}+c^{-1/4})^3(c^{3/4}+1)}{c^{-1}+1}, R(-1/4)=\dfrac{(c^{3/4}+c^{-1/4})(c^{3/4}+1)^3}{c^{-1}+1}\\
L(-1/4)^4&=\dfrac{(c^{3}+c^{-1})^3(c^{3}+1)}{c^{-4}+1}=\dfrac{(c^{4}+1)^3(c^{4}+c)}{c^{4}+1}=(c+1)^8(c^{4}+c),\\
R(-1/4)^4&=\dfrac{(c^{3}+c^{-1})(c^{3}+1)^3}{c^{-4}+1}=\dfrac{(c^{4}+1)(c^{4}+c)^3}{c^{4}+1}=c^3(c^3+1)^3.
\end{split}
\end{equation*}
The third formula follows from $ K(L(-1/4)) =K(R(-1/4))$ of Theorem \ref{thm4-b2cn} and Lemma \ref{lemma1}. The fourth formula arises from the case $ K(L(1/4)) =K(R(1/4))$ and the fifth formula from $ K(L(-1/8)) =K(R(-1/8))$, which the proofs are similar to that of the third formula and omitted.
\end{proof}\qed

\begin{proof}[of Hollmann-Xiang Formula]

From  Helleseth-Zinoviev Formula I, we obtain $ K(c^3(c^3+1)^3)=K(c^9(c^3+1))=K(c^8(c^4+c)) $. From the third formula of Corollary \ref{cor1-from-thm-b2cn}, we have $ K((c+1)^8(c^4+c))=K(c^3(c^3+1)^3)=K(c^8(c^4+c)) $.
\end{proof}\qed

Remark that  Shin-Kumar-Helleseth Formula is the specific case $ K(L(-2)) = K(R(-2))$ of Theorem \ref{thm4-b2cn}, which we omit the proof. 

\subsection{New Kloosterman sum identities and the generalization of  Shin-Kumar-Helleseth Formula}

In this subsection, we establish several identities for Kloosterman sums which generalize  Shin-Kumar-Helleseth Formula.

\begin{thm5}\label{thm5-cbnbk}
Let $ b\in F $, and $ n, k $ be rational numbers such that $ b^n, b^k $ are elements of $ F $. If $ b^{4n}+1\ne 0 $, then, 
\begin{equation*} 
K\biggl(\dfrac{(b^{n+2}+b^k+1)(b^{n+2}+b^n+b^k)^3}{b^{4n}+1}\biggr)=
K\biggl(\dfrac{(b^{n+2}+b^k+1)^3(b^{n+2}+b^n+b^k)}{b^{4n}+1}\biggr).
\end{equation*}
\end{thm5}

\begin{proof}
Let $ k_{1}=b^2+c+1=cb^n+b^k $, then $ c=(b^k+b^2+1)/(b^n+1) $. After an elementary algebraic computation over $ \mathbb{F}_{2} $, we get 
\begin{equation*}
\begin{split}
k_{1}^2&=(b^{n+2}+b^n+b^k)^2/(b^{2n}+1),\\
k_{2}^2&=(b^{n+2}+b^k+1)(b^{n+2}+b^n+b^k)/(b^{2n}+1),\\
(k_{1}k_{2})^2&=(b^{n+2}+b^k+1)(b^{n+2}+b^n+b^k)^3/(b^{4n}+1),\\
(k_{1}k_{2})^2+k_{2}^2&=(b^{n+2}+b^k+1)^3(b^{n+2}+b^n+b^k)/(b^{4n}+1).
\end{split}
\end{equation*}
The actual theorem follows from Lemma \ref{lemma1} and Theorem \ref{thm2}.
\end{proof}\qed

We obtain the following corollaries by affecting particular values to $ n, k $ in Theorem \ref{thm5-cbnbk}, which the proofs are similar to Corollary \ref{cor1-from-thm-b2cn} and omitted.

\begin{cor2}\label{cor2-cbnbk-n1k0}
Set $ n=1, k=0 $  in Theorem \ref{thm5-cbnbk}. Then
\begin{equation*} 
K\biggl(\dfrac{b^3(b^3+b+1)^3}{(1+b)^4}\biggr)=K\biggl(\dfrac{b^9(b^3+b+1)}{(1+b)^4}\biggr).
\end{equation*}
\end{cor2}

\begin{cor3}\label{cor3-cbnbk-n1k2}
Set $ n=1, k=2 $  in Theorem \ref{thm5-cbnbk}. Then
\begin{equation*} 
K\biggl(\dfrac{(b^3+b^2+1)(b^3+b^2+b)^3}{(1+b)^4}\biggr)=K\biggl(\dfrac{(b^3+b^2+1)^3(b^3+b^2+b)}{(1+b)^4}\biggr).
\end{equation*}
\end{cor3}

\begin{cor4}\label{cor4-cbnbk-nm1}
Set $ n=-1$ in  Theorem \ref{thm5-cbnbk}. Then
\begin{equation*} 
K\biggl(\dfrac{(b^{k+1}+b^2+b)(b^{k+1}+b^2+1)^3}{(1+b)^4}\biggr)=K\biggl(\dfrac{(b^{k+1}+b^2+b)^3(b^{k+1}+b^2+1)}{(1+b)^4}\biggr).
\end{equation*}
\end{cor4}

We can obtain analogous results as Theorem \ref{thm4-b2cn} and Theorem \ref{thm5-cbnbk} if we affect one of remainder values from (\ref{eq_KS_k1k2}) to $ k_{1} $ or $ k_{2} $. For instance, set $ k_{2}:=b^2+b+c+\sqrt{c}=b^n\sqrt{c}+c $, we obtain 

\begin{thm6}\label{thm6-bn-csqrt-c}
Let $ b\in F $, and $ n$ be rational number such that $ b^n $ is an element of $ F $. If $ b^{4n}+1\ne 0 $, then, 
\begin{equation*} 
K\biggl(\dfrac{(b^{n+1}+b^2)(b^{n+1}+b^n+b^2+1)^3}{b^{4n}+1}\biggr)=
K\biggl(\dfrac{(b^{n+1}+b^2)^3(b^{n+1}+b^n+b^2+1)}{b^{4n}+1}\biggr).
\end{equation*}
\end{thm6}

\begin{proof}
The proof for the actual theorem is similar to that of Theorem \ref{thm4-b2cn} and Theorem \ref{thm5-cbnbk} and omitted.
\end{proof}\qed

Set $ n=3 $ in the formula of Theorem \ref{thm6-bn-csqrt-c}, we obtain an interesting corollary which the proof is similar to Corollary \ref{cor1-from-thm-b2cn} and omitted:

\begin{cor5}\label{cor5-bn-csqrt-c-n3}
Let $ b\in F $ such that $ b^8+b^4+1\ne 0 $. Then,
\begin{equation*}
K\biggl(\dfrac{(b^3+b^2)(b^3+b+1)^3}{b^8+b^4+1}\biggr)=K\biggl(\dfrac{(b^3+b^2)^3(b^3+b+1)}{b^8+b^4+1}\biggr).
\end{equation*}
\end{cor5}

%

\section{Conclusion}

In this note, we obtained a series of new Kloosterman sum identities  and rediscovered several previously found formulas  with much simpler proof from Theorem 4 of \cite{tor-victor}. The exception is that the formula, $ K(a^5(1+a))= K(a(1+a)^5)$, we believe,  could not be deduced from Theorem 4 of \cite{tor-victor} and Theorem \ref{thm2}.



\end{document}